\documentclass{article}
\usepackage{amsfonts}
\usepackage{amsmath}

\newtheorem{theorem}{Theorem}

\newtheorem{criterion}[theorem]{Criterion}
\newtheorem{definition}[theorem]{Definition}
\newtheorem{example}[theorem]{Example}

\newtheorem{proposition}[theorem]{Proposition}

\newenvironment{proof}[1][Proof]{\noindent\textbf{#1.} }{\ \rule{0.5em}{0.5em}}
\input{tcilatex}

\begin{document}

\title{A U\textsc{nified} A\textsc{pproach} \textsc{of} P\textsc{arameter} E%
\textsc{stimation}}
\author{A\textsc{hmed} G\textsc{uellil}$^{1}$ \textsc{and} T\textsc{ewfik} K%
\textsc{ernane}$^{2}$ \\
$^{1}$\textit{Department of Probability and Statistics, Faculty of
Mathematics }\\
\textit{\ University of Sciences and Technology USTHB,}\\
\textit{BP 32 El-Alia, Algeria}\\
$^{2}$\textit{Department of Mathematics, Faculty of Science}\\
\textit{King Khaled University, Abha Kingdom of Saudi Arabia}\\
e-mail: guellilamed@yahoo.fr, \ tkernane@gmail.com}
\date{}
\maketitle

\begin{abstract}
We introduce a new distance and we use it to parameter estimation purposes.
We observe how it operates and we use in its place the usual methods of
estimation which we call the methods of the new approach. We realize that we
obtain a discretization of the continuous case. Moreover, when it is
necessary to consider truncated data nothing is changed in computations.
\end{abstract}


\section{Introduction}

In the traditional approach of estimation there are three following basic
elements: a family of theoretical probability distributions, an empirical
law and some estimation methods. We choose a method according to its
properties and the problem at hand. The empirical distribution and the
family of theoretical laws are datum of the problem whatever the method
chosen. We propose a new viewpoint where the empirical law corresponding to
a given theoretical one is perceived as being an empirical conditional
distribution with the knowledge of the data. It becomes then an estimate of
the conditional theoretical law knowing the observations before being an
estimation for the theoretical distribution from which it emanated.

We introduce a new distance and we use it to estimate. We observe then how
it operates and use in its place the usual methods of estimation which we
call the methods of the new approach. We notice then that this leads to a
unification of the methods of estimation since we do not make any more
distinction between fixed type-I censored data and complete samples and
between discrete and continuous cases. We thus obtain a considerable
lightening in the procedures of computation in estimation problems. The
distinction in the traditional approach between truncated or type-I censored
data and complete samples is not really justified since all samples are in
fact truncated. Indeed, a sample is not truncated if it covers the totality
of the support of the distribution from which it was drawn, if not it is
truncated. Moreover it is natural to consider that the sample describes only
the parts of the distribution which capture the data. The other parts are
obtained by deduction. Also, the discretization for the continuous case
obtained with the new approach is justified. Indeed, practically all usual
distributions can be reconstituted exactly starting from two or three points
of their graphs. We can then estimate them starting from two or three points
which represent their graphs empirically. In addition to the unification of
several methods of estimation we note that the estimations with the new
measure have the following specific properties. It does not require that the
family of candidate theoretical distributions to be made up of the same type
of laws. There is always a solution which will be acceptable in general. If
the ratios of the frequencies of an empirical distribution coincide with
those of the theoretical one from which it emanated then, from the first we
can find the second with certainty. If the ratios of the frequencies of the
empirical distribution coincide with those of the theoretical one which it
best fits, then the estimations obtained are optimal in the sense that one
cannot improve them. We checked also on some examples, analytically and
numerically, that when we make tending the ratios of the frequencies of the
empirical distribution towards those of the theoretical one, then the
estimates tend towards the true parameters. This last property implies
convergence of the estimators. We prove the convergence of the estimators
obtained with the new measure for a broad class of usual laws. Moreover,
with the new measure we achieve more flexibility in computation compared to
the method of maximum likelihood.

We can distinguish in this paper three different parts. The first is on the
subject of a new distance, presented in section 2. We can be interested and
study it as a mathematical object without necessarily referring to its
applications in statistics. That is a metric which does not have none
equivalent in the theory of mathematics. We noted some of its remarkable
properties, this promises new prospects. The second relates to the use of
this distance in problems of estimation in statistics. That gives birth to a
new method of estimate, presented in section 3. The study suggested in this
part is not at all exhaustive. But the results obtained are already
interesting and encouraging. The third part relates to a new approach of
estimation. We can look at this new approach separately; this is the
discretization of the methods of the continuous case. By adopting it we
widen the field of application of the usual methods of estimation. It is
presented in section 4. In sections 5 and 6 we gave using examples a
practical illustration of the possibilities of the new method and the new
approach of estimation. In section 7 we showed what the users of statistics
gain immediately in the light of our work in comparison with the traditional
approach. Lastly, in section 8 we gave in short a reminder of the whole of
the results obtained.

\section{A New Distance Between Probability Distributions}

In statistics, we use distances to measure the difference between
probability distributions. Usually these distances are conceived in the same
manner, the differences between distributions are almost always expressed by
using variations in geometric sense between their graphs. We introduce a
distance which operates differently. It is based on relativist properties of
probability measures. But its interest is due especially to the fact that it
is not equivalent to usual distances.

\begin{definition}
Consider two probability measures $P$ and $Q$ defined on the same measurable
space $(\Omega ,\mathcal{F})$, $f$ and $g$ being their respective
probability distributions not necessarily with respect to the same measure
and $E$ an event from this space. We say that $f$ and $g$ have same
variations on $E$, if their restrictions on $E$ define the same probability
measure on $E$ equipped with the sigma algebra trace of $\mathcal{F}$ on $E.$
\end{definition}

\begin{example}
Let $f$ be a density of a probability measure $P$ and $E$ an event such that 
$P(E)>0$. The restriction of $f$ on $E$ and the conditional distribution of $%
f$ with respect to $E$ define the same probability measure on $E$ and
consequently they have the same variations on $E.$
\end{example}

\begin{example}
Let $f$ be a probability distribution and $c$ a positive constant. The
functions $f$ and $g=f+c$ have the same variations in the geometric sense
but they do not have the same variations within the meaning of the above
definition.
\end{example}

\begin{proposition}
Let $f$ and $g$ be two probability distributions defined and positives on a
part $E$ not reduced to only one element. If in any point $(x,y)$ of $%
E\times E$, we have%
\begin{equation}
\frac{f(x)}{f(y)}=\frac{g(x)}{g(y)}  \label{ratiofg}
\end{equation}%
then $f$ and $g$ have same variations on $E$.
\end{proposition}

\begin{proof}
If $E$ is discrete the distribution generated by the restriction of $f$ on $%
E $ is $f_{E}=f/\sum_{x\in E}f(x)$ on $E$ and $f_{E}=0$ otherwise. If $x_{0}$
is in $E$ such that $g(x_{0})\neq 0$ then (\ref{ratiofg}) implies that for
all $x$ in $E$, $f(x)=g(x)f(x_{0})/g(x)$. By replacing $f$ in $f_{E}$, we
find the conditional distribution generated by $g$ on $E.$ We obtain then
the result. In the same way, we obtain the result for probability densities
on $\mathbb{R}$ with respect to the Lebesgue measure on $\mathbb{R}$ when $E$
is a subset of $\mathbb{R}$ with positive probability.
\end{proof}

\begin{definition}
Let $f$ and $g$ be two probability distributions and $E$ an event on which
they are strictly positive. If $E$ is discrete and no reduced to only one
element, we call distance in variations between $f$ and $g$ on $E$ the
quantity%
\begin{equation*}
d_{v}(f,g)_{E}=\sum_{\left( x,y\right) \in E}\left\vert \frac{f(x)}{f(y)}-%
\frac{g(x)}{g(y)}\right\vert .
\end{equation*}%
If $E$ is an interval of $\mathbb{R}$ and, $f$ and $g$ are probability
densities on $\mathbb{R}$, with respect to Lebesgue measure $\mu $ on $%
\mathbb{R}$, we call distance in variations between $f$ and $g$ on $E$, the
quantity%
\begin{equation*}
d_{v}(f,g)_{E}=\iint\limits_{E\times E}\left\vert \frac{f(x)}{f(y)}-\frac{%
g(x)}{g(y)}\right\vert \mu (dx)\mu (dy).
\end{equation*}
\end{definition}

Note that $d_{v}$ possesses the properties of symmetry and triangle
inequality. But in the identity property $d_{v}(f,g)_{E}=0%
\Longleftrightarrow f\equiv g$ on $E,$ the equality between $f$ and $g$ must
be understood in the sense that $f$ and $g$ have the same variations on $E$.

Let $d$ be the distance which measures the difference in two points $x$ and $%
y$ between two functions $f$ and $g$ by the quantity $d\left( f,g\right)
\left( x,y\right) =\left\vert f(x)-g(x)\right\vert +\left\vert
f(y)-g(y)\right\vert .$

\begin{proposition}
We have the following property for the distance $d_{v}:$\newline
$d(f,g)(x,y)=0\Longrightarrow d_{v}(f,g)(x,y)=0,$ the converse is not always
true.
\end{proposition}

\begin{proof}
Follows directly from the definitions of $d$ and $d_{v}.$
\end{proof}

\section{New Method of Estimation}

\subsection{Frequency Tables}

Let $\mathcal{F}$ be a family of probability distributions. If it contains
only one type of distribution we say that it is \textit{homogeneous}
otherwise we say that it is \textit{heterogeneous}. A heterogeneous family
can be made up of several types of discrete and absolutely continuous
distributions. Let us consider $f$ in $\mathcal{F}$ and some values $%
y_{1},...,y_{k}$ from its support. We call theoretical table of frequencies
of $f$\ based on $y_{1},...,y_{k}$ or with support $y_{1},...,y_{k}$ the $k$
couples $\left( y_{1},f_{1}\right) ,\left( y_{2},f_{2}\right) ,...,\left(
y_{k},f_{k}\right) $ where $f_{i}=f(y_{i})/%
\sum_{j=1}^{k}f(y_{j}),i=1,2,...,k.$ We note $\bar{f}$ the distribution
defined by this table. We say that the precedent table completely
characterizes the family $\mathcal{F}$ if and only if there is a bijection
between $\mathcal{F}$ and $\mathcal{\bar{F}}$ $=\left\{ \bar{f},\text{ }f\in 
\mathcal{F}\right\} $. In this case, theoretically, from $\bar{f}$ we can
determine $f.$ $\bar{f}$ will be a representative element of $f$ in $%
\mathcal{\bar{F}}$. We call $\mathcal{\bar{F}}$ the family of auxiliary
distributions based on $y_{1},...,y_{k}$ associated to $\mathcal{F}$. We say
also that the $y_{i},$ $i=1,2,\ldots ,k$ form a basis of observations which
characterizes the family $\mathcal{F}$.

\begin{proposition}
Let us consider two laws of probability $f$ and $g$ belonging to a family of
distributions $\mathcal{F}$ and having the same support $E$. If $F$ is a
basis of observations which characterizes the family $\mathcal{F}$ then $%
d_{v}\left( f,g\right) _{F}=0$ implies that $d_{v}\left( f,g\right) _{E}=0.$
\end{proposition}

\begin{proof}
If $d_{v}\left( f,g\right) _{F}=0$ then $\overline{f}=\overline{g}$ where $%
\overline{f}\ $and $\overline{g}$ are the auxiliary distributions of $f$ and 
$g$ respectively based on $F.$ If in addition $F$ constitutes a basis of
observations characterizing $\mathcal{F}$ then, we deduce that $f=g.$
\end{proof}

It should be noted that none of the usual distances has this property and it
is a key idea to justify the use of the methods of point estimation for
discrete case in the continuous one.

\subsection{Estimation}

Let us consider $k$ couples $\left( y_{1},f_{1}\right) ,...,\left(
y_{k},f_{k}\right) $ of a table of empirical frequencies obtained after
grouping the observations of a probability law belonging to a family of
distributions $\mathcal{F}$, with $f_{1}+f_{2}+...+f_{k}=1.$ It will be said
that it empirically characterizes the family $\mathcal{F}$ if the
theoretical frequency table based on the $y_{i}$, $i=1,2,...,k$
characterizes it too. In the sequel our starting point will be always, in
the continuous as in the discrete cases, a table of empirical frequencies,
based on $k$ values $y_{1},...,y_{k},$ constituting a basis of observations
which completely characterizes the studied family. We suppose that it is a
datum of the problem and thus one does not discuss the way of obtaining it,
in particular in the continuous case. We can use for example procedures to
select the optimal number of bins for a regular histogram (see for example
Birg\'{e} and Rozenholc \cite{Birg}). When we use the maximum likelihood
procedure, theoretically nothing prohibits to estimate $n$\ parameters from
a table of empirical frequencies, based on $k$ values where $k$\ is lower or
equal to $n$. But in practice we encounter sometimes difficulties which we
do not expect. In certain cases we note that the results obtained are
completely aberrant. We quote from the literature some paradoxes attached to
the use of the maximum likelihood procedure in these cases (\cite{Josh}).
When we use tables of empirical frequencies whose basis characterizes the
family of theoretical probability distributions which contains the
distribution which we seek we avoid in advance these difficulties\textbf{. }%
We will indicate by $\hat{f}$ the discrete empirical distribution
represented by this table. We notice that it is completely given if the
ratios $f_{i}/f_{j}=\hat{f}(y_{i})/\hat{f}(y_{j})$ $i,j=1,2,...,k$ are known
and if $\hat{f}$ arises from a sample of a given theoretical distribution $f$%
, then from the law of large numbers $\hat{f}(y_{i})/\hat{f}(y_{j})$ tends
to $f(y_{i})/f(y_{j})$ when the sample size tends to infinity. This result
remains valid even when the support $S$ represents a fixed type-I censored
sample. When grouping in classes if one withdraws several classes and their
frequencies, the frequencies of the remaining classes keep this property.
Whether the sample considered is truncated or not and that the distribution
from which it belongs is discrete or absolutely continuous, we can measure
the difference in variations between $\hat{f}$ and a theoretical
distribution $f$ in $y_{1},...,y_{k}$ by%
\begin{equation*}
d_{v}(\widehat{f},f)\left( y_{1},...,y_{k}\right) =\sum_{i,j\in \left\{
1,...,k\right\} }\left\vert \frac{\widehat{f}_{i}}{\widehat{f}_{j}}-\frac{%
f(y_{i})}{f(y_{j})}\right\vert .
\end{equation*}%
Since\textbf{\ }$\widehat{f}$\textbf{\ }converges in probability towards%
\textbf{\ }$f$\textbf{\ }then $d_{v}(\widehat{f},f)$ converges in
probability towards $0$.

Let us consider two probability distribution $f$ and $g$ which does not
belong necessarily to the same type of laws and not equal to zero in $%
y_{1},...,y_{k}$ If $d_{v}(\widehat{f},f)\left( y_{1},...,y_{k}\right)
<d_{v}(\widehat{f},g)\left( y_{1},...,y_{k}\right) ,$ we say that $\widehat{f%
}$ is more close to $f$ than to $g,$ in the sense of $dv$. We thus define a
new method of estimation.

\begin{example}
We simulated $10000$ samples of size $100$ from a binomial distribution $%
\mathcal{B}(8,0.1)$ and $10000$ others from a $\mathcal{B}(15,0.15)$. For
each sample obtained we kept only the observations belonging to $\{0,1,2,3\}$
with their frequencies. Then, starting from the empirical distribution thus
defined we tried to identify the law simulated among the two binomial
distributions considered. The correct distribution is selected for $98,8\%$
of cases when we used samples from the former and for $99,43\%$ of cases
when from the latter.
\end{example}

\begin{example}
We simulated $10000$ samples of size $1000$ from $\mathcal{W}(1.2,1.5)$ and
we omitted the observations below the threshold $1.25$. Each truncated
sample was summarized into $11$ classes. We selected between $\mathcal{W}%
(1.2,1.5)$ and the Gamma distribution $G\left( 2,0.5\right) $ using the
metric $d_{v}.$ The distance $d_{v}$ has selected the correct distribution,
that is $\mathcal{W}(1.2,1.5),$ $98.16\%.$
\end{example}

Let us consider in a problem of estimation, a family of the theoretical laws 
$\mathcal{F}$ and an empirical distribution $\widehat{f}$ with support $%
y_{1},...,y_{k}$ which constitutes a basis of observations characterizing $%
\mathcal{F}$. If it exists $f$ belonging to $\mathcal{F}$ such as $d_{v}(%
\widehat{f},f)\left( y_{1},...,y_{k}\right) =0,$ we say that $f$ is an exact
solution.

\begin{proposition}
The exact solution, when it exists, is optimal in the sense that we cannot
improve it.
\end{proposition}

\begin{proof}
Indeed, in this case there is in $\mathcal{F}$ a distribution whose table of
frequencies coincides exactly with that of $\hat{f},$ it is unique and it is 
$f$.
\end{proof}

\begin{criterion}[of quality]
Let $\hat{f}$ be an empirical distribution and $f$ the theoretical one which
best fits when we estimates by a given method. If $d_{v}(\hat{f},f)=0$ then
according to the preceding proposition the estimate obtained is optimal in
the sense that it is the best possible improvement of the estimation.
\end{criterion}

We have there a quality criterion when it holds, not only it supplants all
the usual criteria but more since it gives a total and definitive guarantee
of the optimality of the estimates. One will further show with examples that
in some cases we can very easily find estimates possessing this property. We
will also show by using examples that, when one makes tending $d_{v}(\hat{f}%
,f)$ towards $0$ the differences between the estimates and the estimated
values tend towards $0$ and at end one obtains their exact values. The
latter property which remains to be proved in the general case implies
immediately convergence of estimates. For the moment there is already the
following result.

\subsection{Convergence in Probability of the Minimum Distance Estimator}

Let $X_{1},...,X_{n}$ a sample with $X_{i}\sim f(x,\theta ),$ $\theta
=\left( \theta _{1},...,\theta _{s}\right) ^{t}\in \Theta \subseteq \mathbb{R%
}^{s},$ with%
\begin{equation}
f(x,\theta )=K(x)\times \exp \left\{ \sum_{k=1}^{s}\theta
_{k}T_{k}(x)+A(\theta )\right\} ,  \label{fam1}
\end{equation}%
$x\in \mathcal{X}\subseteq \mathbb{R},$ where $\mathcal{X}$ is a Borel set
of $\mathbb{R}$ such that $\mathcal{X=}\left\{ x:f(x,\theta )>0\right\} $
for all $\theta \in \Theta .$

The family (\ref{fam1}) is a large family of distributions, one finds there,
for example, the family of the normal laws, and the family of the laws of
Poisson. We assume that the support $\mathcal{X}$ does not depend on $\theta
.$ Denote by $\widetilde{\theta }_{n}$ the estimator by the minimum of
metric $d_{v}$ between the empirical and theoretical distributions $\widehat{%
f}_{n}$ (based on a sample of size $n$) and $f(\cdot ,\theta ),$ that is%
\begin{equation*}
\widetilde{\theta }_{n}=\arg \min_{\theta }d_{v}(f(\cdot ,\theta ),\widehat{f%
}_{n}).
\end{equation*}%
This estimator falls into the class of M-estimators. Using well known
theorems on the convergence of M-estimators (see for example Amemiya \cite%
{Ame}) we will prove that $\widetilde{\theta }_{n}$ converges in probability
to the true parameter.

\begin{proposition}
Let $X_{1},...,X_{n}$ be a sample from the family of distributions (\ref%
{fam1}). If the set of natural parameters $\Theta $ is convex and the true
parameter $\theta $ is an interior point of $\Theta ,$ then the estimator $%
\widetilde{\theta }_{n}$ by the minimum of the distance of variations $d_{v}$
converges in probability to the true parameter $\theta ,$ i.e.,%
\begin{equation*}
\widetilde{\theta }_{n}\overset{P}{\longrightarrow }\theta .
\end{equation*}
\end{proposition}

\begin{proof}
Since we search for a minimum of the criterion function $d_{v},$ it suffices
to show, under the assumptions of the family (\ref{fam1}) and the convexity
of the set $\Theta ,$ that $d_{v}(\theta ,\underline{x})$ seen as a function
of $\theta $ is a convex function (see Amemiya \cite{Ame}). Hence, this
reduces the problem to the convexity of%
\begin{equation*}
\delta _{ij}(\theta )=\left\vert \frac{f(y_{i},\theta )}{f(y_{j},\theta )}-%
\frac{\widehat{f}(y_{i})}{\widehat{f}(y_{j})}\right\vert .
\end{equation*}%
For $\lambda ,\mu \in \mathbb{R}$ with $\lambda +\mu =1$, and $\theta
^{(1)},\theta ^{(2)}\in \Theta ,$ we have%
\begin{equation}
\delta _{ij}(\lambda \theta ^{(1)}+\mu \theta ^{(2)})=\left\vert C_{ij}\exp
\left\{ \sum_{k=1}^{s}\left[ \lambda \theta _{k}^{(1)}+\mu \theta _{k}^{(2)}%
\right] \left( T_{k}(y_{i})-T_{k}(y_{j})\right) \right\} -A_{ij}\right\vert
\end{equation}%
where $C_{ij}=K(y_{i})/K(y_{j})$ and assume that $C_{ij}>0$ and $A_{ij}=%
\widehat{f}(y_{i})/\widehat{f}(y_{j}).$\newline
we have from the convexity of the exponential function that%
\begin{eqnarray*}
\exp \left\{ \sum_{k=1}^{s}\left[ \lambda \theta _{k}^{(1)}+\mu \theta
_{k}^{(2)}\right] \left( T_{k}(y_{i})-T_{k}(y_{j})\right) \right\} &\leq
&\lambda \exp \left\{ \sum_{k=1}^{s}\theta _{k}^{(1)}\left(
T_{k}(y_{i})-T_{k}(y_{j})\right) \right\} \\
&&+\mu \exp \left\{ \sum_{k=1}^{s}\theta _{k}^{(2)}\left(
T_{k}(y_{i})-T_{k}(y_{j})\right) \right\} ,
\end{eqnarray*}%
then%
\begin{equation*}
C_{ij}\exp \left\{ \sum_{k=1}^{s}\left[ \lambda \theta _{k}^{(1)}+\mu \theta
_{k}^{(2)}\right] \left( T_{k}(y_{i})-T_{k}(y_{j})\right) \right\}
-A_{ij}\leq
\end{equation*}%
\begin{equation*}
\lambda C_{ij}\exp \left\{ \sum_{k=1}^{s}\theta _{k}^{(1)}\left(
T_{k}(y_{i})-T_{k}(y_{j})\right) \right\} +\mu C_{ij}\exp \left\{
\sum_{k=1}^{s}\theta _{k}^{(2)}\left( T_{k}(y_{i})-T_{k}(y_{j})\right)
\right\}
\end{equation*}%
\begin{equation*}
-\left( \lambda +\mu \right) A_{ij}\leq \lambda \left[ C_{ij}\exp \left\{
\sum_{k=1}^{s}\theta _{k}^{(1)}\left( T_{k}(y_{i})-T_{k}(y_{j})\right)
\right\} -A_{ij}\right] +
\end{equation*}%
\begin{equation*}
\mu \left[ C_{ij}\exp \left\{ \sum_{k=1}^{s}\theta _{k}^{(2)}\left(
T_{k}(y_{i})-T_{k}(y_{j})\right) \right\} -A_{ij}\right] .
\end{equation*}%
Introducing the absolute value we get%
\begin{equation*}
\delta _{ij}(\lambda \theta ^{(1)}+\mu \theta ^{(2)})=\left\vert C_{ij}\exp
\left\{ \sum_{k=1}^{s}\left[ \lambda \theta _{k}^{(1)}+\mu \theta _{k}^{(2)}%
\right] \left( T_{k}(y_{i})-T_{k}(y_{j})\right) \right\} -\left( \lambda
+\mu \right) A_{ij}\right\vert
\end{equation*}%
\begin{equation*}
\leq \lambda \left\vert C_{ij}\exp \left\{ \sum_{k=1}^{s}\theta
_{k}^{(1)}\left( T_{k}(y_{i})-T_{k}(y_{j})\right) \right\} -A_{ij}\right\vert
\end{equation*}%
\begin{equation*}
+\mu \left\vert C_{ij}\exp \left\{ \sum_{k=1}^{s}\theta _{k}^{(2)}\left(
T_{k}(y_{i})-T_{k}(y_{j})\right) \right\} -A_{ij}\right\vert =\lambda \delta
_{ij}(\theta ^{(1)})+\mu \delta _{ij}(\theta ^{(2)}).
\end{equation*}%
Hence $\delta _{ij}(\theta )$ is a convex function of $\theta ,$ which
implies the convexity of $d_{v}(\theta ,\underline{x})$ seen as a function
of $\theta $ and then the convergence in probability of the minimum of
distance $d_{v}$ estimator.
\end{proof}

\section{New Approach of Estimation}

\subsection{Foundation}

Let us consider in a problem of estimation the family of theoretical
distributions $\mathcal{F}$ and an element $f$\ belonging to $\mathcal{F}$.
We have in an obvious way, $d_{v}(\widehat{f},f)\left(
y_{1},...,y_{k}\right) =d_{v}(\widehat{f},\bar{f})\left(
y_{1},...,y_{k}\right) $ where $\bar{f}$ is the representative of $f$ in $%
\mathcal{\bar{F}}$, $\mathcal{\bar{F}}$ being the family of auxiliary
distributions based on $y_{1},...,y_{k}$, associated to $\mathcal{F}.$ $\bar{%
f}$ is a discrete probability distribution with same support as $\hat{f}$
and depend on the same parameters of $f.$ If the theoretical table of
frequencies based on $y_{1},...,y_{k}$ characterizes completely the family $%
\mathcal{F}$ then the determination of $f$ is equivalent to the
determination of $\bar{f}$. When $\mathcal{F}$ is homogeneous, for
determining $\bar{f}$, instead of $dv$ we can also make use of the usual
methods (method of moments, method of maximum likelihood, Bayesian Methods,
... etc.). Then they will be called the methods of the new approach. When
proceeding in this way, all occurs as if one replaces the family of the
theoretical distributions $\mathcal{F}$ by the corresponding family $%
\mathcal{\bar{F}}$. We note also what follows:

\textbf{1.} In discrete case, if the usual methods of estimation are used it
is as if one estimates in a traditional way starting from truncated samples.
This supposes that it is considered that any sample which does not
completely cover the support of the distribution from which it is resulting
is truncated in a deterministic way, the truncation being the parts which do
not appear in the observations.

\textbf{2.} In continuous case, often in practice one associates with the
sample of observations an optimal discrete distribution in a certain way and
one uses it to estimate. Then when replacing $dv$ by the usual methods we
obtain a discretization of the continuous case.

\textbf{3.} In discrete case $\bar{f}$ represents the conditional
distribution of $f$ knowing the observations $y_{1},...,y_{k}$. In the
continuous case $\bar{f}$ is calculated in a similar manner. It seems that
there also it has the same interpretation except that this type of
calculation does not exist in the theory of probability.

For reason of coherence only with what has just been said in 1, 2 and 3, we
propose to view the empirical distribution as being the conditional
empirical distribution knowing the observations, since it is calculated
knowing the observations, even if that is not obvious in the continuous
case. One then conceives it more easily as being an estimate of $\bar{f}$
before being for $f$.

\section{Analytical computation}

In this part we will organize a discussion around some very simple examples
to try to reveal the specificity of the new approach and its contribution
compared to the traditional one. Let us consider a table of frequencies
based on two observations $x$ and $y$ with their respective frequencies $%
n_{1}$ and $n_{2}$. Starting from such table, with the new method one can
estimate only one parameter. Such table characterizes practically all the
families of usual laws when one has to estimate only one parameter. We can
obtain such a table when the sample considered is not truncated but of small
size or is truncated and grouped in two classes only. In the light of the
new distance we will see in the examples which follow that, according to
whether one estimates only one parameter or two simultaneously, even if the
sample is not of small size, it will be henceforth preferable to group it in
two or three classes only because one can gain in the precision of the
estimations. Indeed, the two or three points obtained have more weight to
represent the theoretical points of the distribution which they describe
empirically and the method of estimation with $d_{v}$ practically always
gives in this case an optimal solution in the most general meaning.

\subsection{Estimation of the parameter of the exponential distribution}

Assume we want to estimate, from the preceding table, the probability
density $f_{\lambda }$ given by $f_{\lambda }$ $(x)=\lambda e^{-\lambda x}$
if $x>0$ and $f_{\lambda }(x)=0$ otherwise, $\lambda >0$, and $F$ denotes
its cdf.

\textbf{a.} Suppose it is a summary of a not truncated sample. Then the
estimators of $\lambda $ by the methods of maximum likelihood of the
classical approach $\hat{\lambda}$ and the new one $\hat{\lambda}_{N}$ are
respectively: $\hat{\lambda}=\left( n_{1}+n_{2}\right) /\left(
n_{1}x+n_{2}y\right) $ and $\hat{\lambda}_{N}=\left( \log \left(
n_{1}\right) -\log (n_{2})\right) /\left( y-x\right) .$As we can see, in
general $\hat{\lambda}$ is different from $\hat{\lambda}_{N}$. When we
compute $\tilde{\lambda},$ the estimation obtained using $d_{v},$ we find
that it is equal to $\hat{\lambda}_{N}.$ $\tilde{\lambda}$ is here optimal
in the general sense. If

\begin{equation*}
\frac{n_{1}}{n_{2}}=\frac{f(x)}{f(y)}+\varepsilon
\end{equation*}%
then 
\begin{equation*}
\hat{\lambda}_{N}\left( \varepsilon \right) =\lambda +\varepsilon k
\end{equation*}%
$k$ being a constant. Then $\hat{\lambda}_{N}\left( \varepsilon \right) $
tends towards $\lambda $ when $\varepsilon $ tends towards $0.$ We can check
that the difference between $\lambda $ and $\hat{\lambda}\left( \varepsilon
\right) $ does not tend towards $0$ when $\varepsilon $ tends towards $0$.
If the sample size tends towards infinity then, from the law of large
numbers, the differences between the ratios of the empirical relative
frequencies and those theoretical which their correspond tend towards $0$
and consequently $\hat{\lambda}_{N}$ tends to $\lambda .$ But one can have
these variations close to $0$ same for samples of finite sizes. It is
noticed that the first solution here is always acceptable but the second
not. The second is not acceptable only if there are anomalies in the sample
of observations and then one is warned. We are not able to detect the sample
deficiency from the first. The second is not acceptable when $x<y$ and $%
n_{2}<n_{1}$ or conversely, but it is not what one expects, since the
exponential law being decreasing, $x<y$ we must have $n_{1}>n_{2}$. Now if
in a problem the preceding exact solution is not acceptable and we have to
propose an estimate of $\lambda ,$ that is always possible with the new
method. Put

\begin{equation*}
\alpha \left( \lambda \right) =\left\vert \frac{f(x)}{f(y)}-\frac{n_{1}}{%
n_{2}}\right\vert +\left\vert \frac{f(y)}{f(x)}-\frac{n_{2}}{n_{1}}%
\right\vert \text{ and }E=\left\{ \alpha \left( \lambda \right) ,\lambda
>0\right\}
\end{equation*}%
$E$ is a part of $\mathbb{R}$ which is bounded below by $0.$ It admits then
a lower bound say $\alpha _{0}$. If $\alpha _{0}$ is in $E$ then there is $%
\lambda _{0}>0$ such that $\alpha (\lambda _{0})=\alpha _{0}.$ In this case
the estimation of $\lambda $ is $\lambda _{0}.$ If $\alpha _{0}$ is not in $%
E $ then, whatever the strictly positive integer $n,$ there exists $\lambda
>0$ such that $\left\vert \alpha (\lambda )-\alpha _{0}\right\vert <1/n.$
Put $A_{n}=\left\{ \lambda >0/\text{ }\left\vert \alpha (\lambda )-\alpha
_{0}\right\vert <1/n\right\} .$ $A_{n}$ is a decreasing sequence and then
there exists $A_{0}$ such that $\lim_{n\rightarrow \infty }A_{n}=A_{0}.$ In
this case, each value $\lambda $\ from $A_{0}$ can be considered as an
estimation of $\lambda $ with the new approach.

\textbf{b.} Assume now that the table given is that of a fixed type-I
censored data. For example in a not truncated grouped data one kept only the
centers of two classes and their corresponding frequencies. With the new
approach the table is enough and the solution is exactly the same as
previously. But in this case the preceding estimate of the traditional
approach is not valid here. One must use the methods of truncated data. One
then needs the part of the support of $f$ represented here by $x$ and $y.$
To be able to carry out calculations let us suppose that this table is the
summary of the observations falling into the interval $[0,c]$ with $c>0$.
That is a right truncated sample. We consider the observed likelihood%
\begin{equation*}
L_{obs}=\left( \frac{f(x)}{F(c)}\right) ^{n_{1}}\left( \frac{f(y)}{F(c)}%
\right) ^{n_{2}}.
\end{equation*}%
We have to consider that $n_{T}$ observations are greater than $c$ and have
been discarded, but $n_{T}$ is unknown. In order to compute the complete
likelihood we have to determine the conditional distribution of $n_{T}$
given that the observations follows an exponential distribution to be able
to implement the EM\ algorithm which require the computation of the
conditional expectation of the complete log-likelihood function. It is then
not possible to have an analytic solution and a recursive procedure is used
to achieve a numerical solution. In general it is not always easy to use the
method of maximum likelihood as let it believe the examples on the usual
laws. Although Maximum likelihood estimators have good statistical
properties in large samples, they often cannot be reduced to simple
formulas, so estimates must be calculated using numerical methods.

\subsection{Estimation of the parameters of a normal distribution}

Let us consider a normal law $N\left( m,\sigma \right) .$

\subsubsection{Estimation of the average}

Solving the following equation in $m:$

\begin{equation*}
\frac{n_{1}}{n_{2}}-\frac{\exp \left( -\frac{\left( x-m\right) ^{2}}{2\sigma
^{2}}\right) }{\exp \left( -\frac{\left( y-m\right) ^{2}}{2\sigma ^{2}}%
\right) }=0
\end{equation*}

we obtain

\begin{equation*}
\tilde{m}=\frac{1}{-\frac{x}{\sigma ^{2}}+\frac{y}{\sigma ^{2}}}\left( -\ln 
\frac{n_{1}}{n_{2}}-\frac{1}{2}\frac{x^{2}}{\sigma ^{2}}+\frac{1}{2}\frac{%
y^{2}}{\sigma ^{2}}\right)
\end{equation*}%
It should be noted that $\tilde{m}$ is function of $\sigma .$ When solving
precedent equation after replacing $\left( n_{1}/n_{2}\right) $ by $\left(
f(x)/f(y)\right) +\varepsilon $, we obtain:

\begin{equation*}
\tilde{m}\left( \varepsilon \right) =\frac{1}{-\frac{x}{\sigma ^{2}}+\frac{y%
}{\sigma ^{2}}}\left( -\ln \left( \frac{e^{-\frac{\left( x-m\right) ^{2}}{%
2\sigma ^{2}}}}{e^{-\frac{\left( y-m\right) ^{2}}{2\sigma ^{2}}}}%
+\varepsilon \right) -\frac{1}{2}\frac{x^{2}}{\sigma ^{2}}+\frac{1}{2}\frac{%
y^{2}}{\sigma ^{2}}\right)
\end{equation*}%
where $\lim\limits_{\varepsilon \rightarrow 0}\tilde{m}\left( \varepsilon
\right) =m.$

\subsubsection{Estimation of the Variance}

Solving the following equation in $\sigma $,

\begin{equation*}
\ln \frac{n_{1}}{n_{2}}=-\frac{\left( x-m\right) ^{2}}{2\sigma ^{2}}+\frac{%
\left( y-m\right) ^{2}}{2\sigma ^{2}}
\end{equation*}

we have:

\begin{enumerate}
\item If $\frac{n_{1}}{n_{2}}=1$ and $-2mx+2my+x^{2}-y^{2}=0,$any value $%
\sigma $ belonging to $\mathbb{R}$ is solution.

\item If $\frac{n_{1}}{n_{2}}=1$ and $-2mx+2my+x^{2}-y^{2}\neq 0,$ there is
no solution.

\item If $\frac{n_{1}}{n_{2}}\neq 1,$one obtains:
\end{enumerate}

\begin{equation*}
\tilde{\sigma}=\left\vert \frac{1}{2\ln \frac{n_{1}}{n_{2}}}\sqrt{2}\sqrt{%
2mx\ln \frac{n_{1}}{n_{2}}-2my\ln \frac{n_{1}}{n_{2}}-x^{2}\ln \frac{n_{1}}{%
n_{2}}+y^{2}\ln \frac{n_{1}}{n_{2}}}\right\vert
\end{equation*}

If $\frac{n_{1}}{n_{2}}=\frac{f(x)}{f(y)}+\varepsilon ,$ one obtains $%
\lim\limits_{\varepsilon \rightarrow 0}\tilde{\sigma}\left( \varepsilon
\right) =\sigma .$

\subsection{Remarks}

\textbf{1. }As shown in the examples above, if there is a table of
frequencies based on two observations and one estimates only one parameter,
then with $dv$ one easily obtains optimal estimates in the most general
sense of the term. It is not always easy when the table is based on $k$
observations $y_{1},...,y_{k}$ with $k\geq 3$. If the table is thus formed
and that we cannot determine a total exact solution one proposes to take the
various couples of possible observations in $\left\{ y_{1},...,y_{k}\right\} 
$ and to determine the exact solution each time when it is possible and
approached otherwise. Each estimation is weighted by the sum of the
frequencies of the elements of the couple and we calculate their mean. For
example in the case of the first example if there are exact solutions for
the various couples we take $\tilde{\lambda}=\left(
1/\sum\limits_{i,j=1,i\neq j}^{k}\left( n_{i}+n_{j}\right) \right)
\sum\limits_{i,j=1,i\neq j}^{k}\left( n_{i}+n_{j}\right) \frac{\left( \ln
\left( n_{i}\right) -\ln (n_{j})\right) }{y_{j}-y_{i}}.$ We notice that here
for each couple the estimation converges towards the true value when the
differences between the ratios of the empirical relative frequencies and
corresponding theoretical ones tend towards $0$, then it is the same for the
latter.

\textbf{2.} In the first example we have obtained the same solution with $dv$
and the method of maximum likelihood of the new approach. It is not an
isolated case. We noted in various examples considered in this document,
when we estimate only one parameter, they always give concordant results.

\section{Numerical Example}

Even in the discrete case the two approaches are different since, contrary
to the traditional one, with the new we do not distinguish truncated samples
from those not truncated. In traditional approach of truncated samples all
parts of the support of the estimated distribution which are supposed to be
observed are used in calculations through the conditional theoretical
distribution. With the new one we use only the observations. Now, if we
consider the samples which do not cover all the support of the distribution
from which they emanated are truncated, the truncations being the parts
which do not appear in the observations and we apply the traditional
approach, we fall in the new one. For this reason we do not insist on the
discrete case, we give only examples concerning the continuous case. It is
not easy to present a comparative study of the numerical results of the two
approaches, since to the same estimate of the new it corresponds two
estimates of the traditional according to whether it is considered that the
sample is truncated or not. In addition, in the traditional approach when
the sample is truncated the nature of truncation is used in calculations.
Then the frequency table, without indication of the parts observed, is not
enough. It is necessary at each time to indicate the intervals represented
by the observations in the table. For all these reasons we present the
estimates of the two approaches only when that makes better to underline the
specificity of the new one. For example, we simulated synthetic data of size 
$400$ from the standard normal distribution and we grouped them into $11$
classes represented by the observations $y_{1},...,y_{11}$ and their
frequencies. We obtain $y_{3}=$ $-1.5331,$ $y_{6}=0.0386$ and $y_{8}=1.0863$
with their respective absolutes frequencies $n_{3}=23,$ $n_{6}=89$ and $%
n_{8}=43.$ In the table presented hereafter, in the part before the line of $%
n_{8}$ we consider the two observations $y_{3}$ and $y_{6}$. The distance $%
dv $ in these two points between the empirical distribution and the standard
normal distribution is null as one takes $n_{3}=27500$ and $n_{6}=89000$. We
fix then $n_{6}=89000$ and give ascending values for $n_{3}$, more and more
near to $27500$ as indicated in the table and we estimate $m$ when $\sigma $
is known and $\sigma $ when $m$ is known. At each time we estimate them with
the method of minimal distance with $dv$, the method of moments of the new
approach and the method of maximum likelihood of the classical approach. We
note estimates obtained with $dv$ and with maximum likelihood of the new
approach respectively by $\widetilde{m}$ and $\hat{m}_{Mnew}$ for average
and $\tilde{\sigma}$ and $\hat{\sigma}_{Mnew}$ for the standard deviation
and we note $\hat{m}_{CLH}$ and $\hat{\sigma}_{CLH\text{ }}$ those obtained
with the classical maximum likelihood procedure for truncated samples. For
this last, the observed part is assumed to be $\left[ -1.7951,-1.2712\right[
\cup \left[ -0.22335,0.30055\right[ .$

$%
\begin{tabular}{|l|c|c|c|c|c|}
\hline\hline
\multicolumn{6}{|c|}{$y_{3}=\mathbf{-1.5331},$ $y_{6}=\mathbf{0.038690}%
,y_{8}=\mathbf{1.0863},n_{6}=\mathbf{89000}$} \\ \hline\hline
$\mathbf{n}_{3}$ & $\mathbf{23000}$ & $\mathbf{24000}$ & $\mathbf{26000}$ & $%
\mathbf{27000}$ & $\mathbf{27500}$ \\ \hline
$\tilde{m}$ & \multicolumn{1}{|l|}{$0.113\,69$} & \multicolumn{1}{|l|}{$%
0.08661$} & \multicolumn{1}{|l|}{$0.03568\,$} & \multicolumn{1}{|l|}{$%
0.01167 $} & \multicolumn{1}{|l|}{$-0.000001$} \\ \hline
$\hat{m}_{Mnew}$ & $0.11369$ & $0.08661$ & $0.03568$ & $0.01167$ & $%
-0.000001 $ \\ \hline
$\hat{m}_{CLH}$ & \textbf{0.110\thinspace 75} & \textbf{0.08444} & \textbf{%
0.03478} & \textbf{0.01128} & \textbf{0.000155} \\ \hline
$\tilde{\sigma}$ & \multicolumn{1}{|l|}{$0.931\,64$} & \multicolumn{1}{|l|}{$%
0.946\,64$} & \multicolumn{1}{|l|}{$0.976\,94$} & \multicolumn{1}{|l|}{$%
0.992\,28$} & \multicolumn{1}{|l|}{$1.0$} \\ \hline
$\hat{\sigma}_{Mnew}$ & \multicolumn{1}{|l|}{$0.931\,64$} & 
\multicolumn{1}{|l|}{$0.94664$} & \multicolumn{1}{|l|}{$0.976\,94$} & 
\multicolumn{1}{|l|}{$0.99228$} & \multicolumn{1}{|l|}{$1.0$} \\ \hline\hline
$\hat{\sigma}_{CLH}$ & \textbf{0.92171} & \textbf{0.93701} & \textbf{0.967796%
} & \textbf{0.98335} & \textbf{0.991165} \\ \hline\hline
$\mathbf{n}_{8}$ & $\mathbf{43000}$ & $\mathbf{44444}$ & $\mathbf{47273}$ & $%
\mathbf{48214}$ & $\mathbf{49371}$ \\ \hline
$\tilde{m}$ & $-0.02224$ & $-0.017549$ & $-0.00785$ & $-0.00762$ & $0.000002$
\\ \hline
$\hat{m}_{Mnew}$ & $0.03676\,3$ & $0.05190\,7$ & $0.088443$ & $0.102\,94$ & $%
0.0000005$ \\ \hline
$\tilde{\sigma}$ & $0.917\,67$ & $0.935\,46$ & $0.97180$ & $0.987\,16$ & $%
1.0 $ \\ \hline
$\hat{\sigma}_{Mnew}$ & $1.068\,9$ & $1.108\,0$ & $1.\,\allowbreak 196\,8$ & 
$1.242$ & $1$ \\ \hline
\end{tabular}%
$

In the part after the line of $n_{8}$ we estimate simultaneously $m$ and $%
\sigma $ by the method of the minimal distance with $dv$ and the method of
moments of the new approach starting from the observations $y_{3}$, $y_{6}$
and $y_{8}$ by fixing the frequency of $n_{8}$ $=89000$ and while taking for 
$n_{3}$ and $n_{6}$, the frequencies indicated. Then we observe what occurs
when we make tending the differences between the ratios of the empirical
frequencies and the corresponding theoretical frequencies towards $0.$ It is
noticed that in the various examples considered, when we estimate only one
parameter, the various methods of the new approach agree completely. But it
is not the case when one estimates simultaneously two parameters. In the
table above, when we estimate simultaneously $m$ and $\sigma $ with the
method of the moments of the new approach or the method of minimal distance
with $dv$, when the ratios of the empirical frequencies coincide exactly
with the corresponding theoretical ones we obtain their exact values. But
with the method of moments, as we can see, the difference between the
estimated parameters and their true values does not decrease necessarily
when the difference between these ratios decreases as with the method of the
minimal distance with $dv$. It seems that this property is specific to the
estimation with $dv$. Here, in the various estimates with $dv$, at each
time, the distance within the meaning of $dv$ between the empirical
distribution considered and the one to which it leads is null. Consequently
the estimates with $dv$ in that table are optimal in the most general
meaning.

\section{Comparison of the two approaches}

A more thorough study is needed to compare the two approaches of estimation
than only one section. But, by putting ourselves in the viewpoint of users
of statistics, we can try to characterize what is achieved with the new
approach at various levels.

\subsection{Procedures}

We place at disposal of statisticians all the usual methods of estimation
and a new one. The remarkable fact with the new approach is that it occurs
as if all is discrete except the need for grouping observations into classes
in the continuous case. moreover, when it is necessary to consider fixed
type-I censoring nothing change in computations. With this unification of
several methods of estimation we obtain a considerable lightening of
procedures compared to the traditional approach.

\subsection{Computations}

With the new approach, since all is discrete, there is no more the usual
difficulties related to the integral calculus. With the method of maximum
likelihood of the traditional approach or the new one, sometimes we
encounter great difficulties when one must estimate several parameters
simultaneously. But with the method of the minimal distance with $d_{v}$ one
can always easily propose an acceptable solution.

\subsection{Credibility of estimates.}

The statistician can now estimate with various methods, those of the
traditional approach and of the new. If he obtains two different appreciable
results it must decide for one of them. Usually we do not decide in this way
since in the traditional approach we do not have criteria which give
guarantees on a given specific evaluation. We have only criteria which give
guarantees on average or asymptotically or by confidence interval. In this
spirit, to make admitting the new approach we should prove that it makes
possible to obtain estimations better relatively to these criteria compared
to those usually obtained. If one places itself in this spirit then, it is
useless to continue because, for example, one cannot find better than the
empirical average to estimate the average of the normal law. Of course
nothing prevents us from also looking at the usual criteria in the new
approach but there are new elements. One can henceforth in certain cases,
without determining the estimator, affirming with certainty that the point
estimation obtained with the new method is better than that obtained with
maximum likelihood procedure. In other cases one can give estimators and
without studying their properties one can affirm that one cannot improve
them. Indeed, when the distance, within the meaning of $d_{v}$, between a
given empirical distribution and the theoretical one which best fits is
null, the estimate obtained is optimal in the general sense. It is noticed
that when the distance within the meaning of $d_{v}$ between a given
empirical distribution and the one we obtain by the method of the minimal
distance with $d_{v}$ is not null, the solution obtained is regarded as
optimal only within the meaning of the $d_{v}$. In this case perhaps it is
optimal in the most general sense what must then be specified. This question
remains to be studied.

\section{Conclusion}

We introduced a new distance and we proposed an new approach of the
estimation.

\textbf{1. The New distance.}

We introduced a new distance and we used it in parameter estimation where we
noticed what follows.

a. One can estimate even when the family of candidate theoretical
distributions is not homogeneous and there is always a solution which will
be acceptable in general.

b. Given a discrete empirical distribution associated to a sample belonging
to a theoretical one,

- If the ratios of frequencies of the first coincide with those of the
second we found exactly the latter.

- If the ratios of the frequencies of the first coincide with those of the
theoretical one which best fits, then the estimations obtained are optimal
in the sense that one cannot improve them.

- We showed on some examples that if we make tending the ratios of the
frequencies of the first towards the corresponding theoretical ones of the
second, then the estimations tend towards the true parameters. This implies
immediately the convergence of the estimators. We showed the convergence in
probability of the estimator for a broad class of usual laws.

c. We introduced a quality criterion, when it holds, it is stronger than of
checking all the usual criteria together and we showed on some examples that
in certain cases we can determine easily estimations which check it.

In addition we note a certain flexibility in calculations with $dv$ compared
to the method of the maximum likelihood.

\textbf{2. The New approach.}

We proposed an new approach of parameter estimation. When it is applied it
works as if all is discrete except the need for grouping the observations in
bins in continuous case. Since all is discrete there is no more the usual
difficulties related to integral calculus. moreover, when it is necessary to
consider fixed type-I censoring nothing is changed in computations. This
unification of several methods of estimation leads to a lightening of the
procedures compared to the traditional approach.

\end{document}